\documentclass[a4paper,12pt]{article}
\usepackage[utf8]{inputenc}
\usepackage[]{amsmath}
\usepackage[]{amssymb}
\usepackage[]{a4wide}
\usepackage[]{amsthm}
\author{Magdalena Zielenkiewicz}
\title{Integration over homogenous spaces for classical Lie groups using iterated residues at infinity}
\date{}
\def\T{\mathbb{T}}
\def\Grass#1#2{\textrm{Grass}_{#1}(\mathbb{C}^{#2})}
\def\VR{V_{\mathcal{R}}}
\def\VQ{V_{\mathcal{Q}}}
\def\R{\mathcal{R}}
\def\lie#1{\mathfrak{#1}}
\newtheorem{thm}{Theorem}
\newtheorem{cor}{Corollary}
\newtheorem{formula}{Formula}
\theoremstyle{definition}
\newtheorem{example}{Example}

\begin{document}
\maketitle
\begin{abstract}
Using the Berline-Vergne integration formula for equivariant cohomology for torus actions, we prove that integrals over Grassmannians (classical, Lagrangian or orthogonal ones) of characteristic classes of the tautological bundle, can be expressed as iterated residues at infinity of some holomorphic functions of several variables. The results obtained for these cases, can be expressed as special cases of one formula, involving the Weyl group action on the characters of the natural representation of the torus.
\end{abstract}

\section{Introduction}
The goal of this paper is to present a way of expressing integrals over Grassmannians with torus action as iterated residues at infinity of holomorphic functions. Localization theorems for compact spaces with torus action provide a useful tool, which allows to express some cohomological invariants using only the fixed point set of the action.
For example, if $X$ is a compact manifold with torus action and the fixed point set of the action is finite, the Berline-Vergne formula expresses the integral over $X$ of an element of equivariant cohomology as a sum of local contributions:
 \[\int_{X} \alpha  = \sum_{p \in X^\T} \frac{i_p^*\alpha}{e_p}.\]
Here, $\alpha \in H^*_\T(X)$ is an element of equivariant cohomology of $X$, $e_p$ is the equivariant Euler class of the tangent bundle and $\T$ denotes the n-dimensional torus. If the space $X$ satisfies additional assumptions (e.g. the odd dimensional cohomology vanishes), then every closed $\T$-invariant differential form can be extended to an element of $H^*_\T(X)$. \cite{ginzburg} \newline

In the case of complex Grassmannian, the Berline-Vergene formula has the form:
\[ \int_{\Grass{m}{n}} \phi(\R) = \sum_{p_\lambda} \frac{\phi_{|_{p_\lambda}}}{e_{p_\lambda}} .\]
Now, $\phi(\R)$ is a characteristic class of the tautological bundle $\R$. To every characteristic class we can associate a symmetric polynomial $V$, in Chern roots of the bundle $\R$.  I prove a following formula for the above integral, that has been claimed in the paper of A.Weber \cite{weber}:
\[ \int_{\Grass{m}{n}} \phi(\R) = \frac{1}{m!} Res_{\mathbf{z} = \infty} \frac{V(z_1,...,z_m)\prod_{i \neq j}(z_i - z_j)}{\prod_{i=1}^n\prod_{j=1}^m(t_i - z_j)}. \]
A similiar result can be obtained for other types of Grassmannians. In the case of Lagrangian Grassmannian $LG(n)$, the expression for an integral, in the form of an iterated residue, is the following:
\[\int_{LG(n)} \phi(\R) =  Res_{\mathbf{z} = \infty} \frac{1}{n!}\frac{V(z_1,...,z_n)\prod_{i \neq j}(z_j - z_i)\prod_{i < j}(z_i+z_j)}{\prod_{i=1}^n\prod_{j=1}^n(t_i - z_j)(t_i + z_j)}\]
The formulas for orthogonal Grassmannians are very similiar to the one above, and are given in chapter \ref{orthogonal}. Berczi and Szenes in their paper \cite{berczi} found a formula for an integral over flag variety (\cite{berczi}, chapter 6.3).
All the obtained formulas (for classical, Lagrangian and orthogonal Grassmannians and the one derieved by Berczi and Szenes) can be written as one equality, involving the action of the Weyl group on the characters of the natural representation of torus action. The result is presented in chapter \ref{general}.

\section{Preliminaries}
Throughout this paper, we will use some well-known facts about equivariant cohomology for torus action, and basic knowledge of theory of analytic functions (we use iterated residues of  holomorphic functions at infinity).

\subsection{Berline-Vergne formula}
Let $\T$ be the n-dimensional torus $\T=(\mathbb{C}^*)^n$, acting on a compact space $X$. The characters of the action can be identified with elements of equivariant cohomology of $X$ - to every character $\chi$ we associate the first equivariant Chern class of the line bundle defined by this character, $c_1(\chi) \in H_\T^2(pt)$. Thus we can consider the group of nonzero characters $\T^{\#}\setminus \{0\}$ as a subset of $H^*_\T(X)$.

\begin{thm}[Atiyah, Bott] Let $X$ be a compact $\T$-space. The inclusion  $i: X^\T \hookrightarrow X$ induces an isomorphism,
\[ i^*: H_\T^*(X)[(\T^{\#} \setminus\{ 0\})^{-1}] \stackrel{\sim}{\longrightarrow} H_\T^*(X^\T)[(\T^{\#} \setminus\{ 0\})^{-1}],\]
after localizing with respect to the multiplicative system consisting of finite products of elements $c_1(\chi)$, for $\chi \in \T^{\#}\setminus \{0\}$. \cite{quillen}
\end{thm}

If $X$ is a compact manifold, and the fixed point set is finite, the Atiyah-Bott localization theorem has the following consequence:

\begin{thm}[Berline, Vergne] Suppose $X$ is a compact manifold and $\# X^\T < \infty$. For $x \in H_\T^*(X)$ we have
\[ \int_{X} x  = \sum_{p \in X^\T} \frac{i_p^*x}{e_p},\]
where $e_p$ is the equivariant Euler class of the tangent bundle at the fixed point $p$, and $i_p^*$ is the restriction of the inclusion $i:X^\T \hookrightarrow X$ to point $p$. \cite{ab}
\end{thm}

The equivariant Euler class at the fixed point $p$ is given by the product of weights of torus action on the tangent space. If $X$ is a homogenous space $G/P$, where P is a parabolic subgroup (like in the case of Grassmannians and flag varieties), the tangent space at $1$ is canonically isomorphic to $\lie{g}/\lie{p}$ and the weights of the torus action are the posivite roots~$\Phi^+ \setminus \Phi_P^+$.

\subsection{Iterated residues at infinity}
Let $\omega_1,...,\omega_k$ be affine forms on $\mathbb{C}^n$, and let $h(z_1,\dots,z_n)$ be a holomorphic function. We define the iterated residue at infinity as
\[ Res_{z_1 = \infty}...Res_{z_n = \infty} \frac{h(\mathbf{z})d\mathbf{z}}{\prod_{i=1}^k \omega_k} : = \big(\frac{1}{2\pi} \big)^n \int_{|z_1| = R_1} \cdots \int_{|z_n| = R_n}  \frac{h(\mathbf{z})d\mathbf{z}}{\prod_{i=1}^k \omega_k}, \]
where $\mathbf{z} : =  (z_1,\dots,z_n) $ and $d\mathbf{z} : =  dz_1\wedge...\wedge dz_n$, and the intergals are taken over sufficently big circles, satisfing $1 \ll R_1 \ll ... \ll R_n$. 

\section{Results}\label{results}
The goal of this paper is to show how integrals over various Grassmannians (classical, Lagrangian and orthogonal) can be expressed as interated residues at infinity of rational functions. The starting point is the Berline-Vergne formula, and the following example:
\begin{example}[Berline-Vergne formula for $\mathbb{CP}^n$]
Consider $\mathbb{CP}^n$ with the standard torus action: \[(t_0,...,t_n) \cdot [x_0:...:x_n] \mapsto [t_0x_0:...:t_nx_n].\] The fixed points are the coordinate lines $p_i = [0:\dots:1:\dots:0]$. The tangent space at point $p_i$ is
\[T_{p_i}\mathbb{CP}^{n} = Hom(p_i, p_i^{\perp})\]
and therefore the Euler class is equal to \[e_{p_i} = \prod_{\substack{l=0 \\ l \neq i}}^n (t_l - t_i).\]
Using Berline-Vergne formula to compute the integral of a characteristic class $\phi(\mathcal{R})$, given by a polynomial $V(c_1(\mathcal{R}))$, we obtain 
\[\int_{\mathbb{CP}^n} \phi(\mathcal{R}) =  \sum_{i=0}^n \frac{V(t_i)}{\prod_{j \neq i}(t_j - t_i)} . \]
One can easily write the righthand-side as a residue:
\[ \sum_{i=0}^n \frac{V(t_i)}{\prod_{j \neq i}(t_j - t_i)} = \sum_{i=0}^n - Res_{z=t_i}\frac{V(z)}{\prod_{j \neq i}(t_j - z)(t_i-z)} = Res_{\mathbf{z}=\infty}\frac{V(z)}{\prod_{j=0}^n(t_j - z)},\]
where the last equality holds by the Residue Theorem. 
\end{example}

Using similar calculations, one can show an analogous formula for integral over the Grassmannian of $m$-dimensional subspaces in $\mathbb{C}^n$. The fixed points of the torus action are the coordinate subspaces 
\[p_{\lambda} = Span\{\epsilon_{\lambda_1},...,\epsilon_{\lambda_m}\},\]
where $\lambda =  (\lambda_1<...<\lambda_m)$ is a partition of $n$. The Euler class is equal to
\[ e_{p_\lambda} = \prod_{k \in \lambda, l \notin \lambda}(t_l - t_k) .\]
Therefore, the Berline-Vergne formula has the form:
\[ \int_{\Grass{m}{n}} \phi(\mathcal{R}) = \sum_{p_\lambda} \frac{\phi_{|_{p_\lambda}}}{e_{p_\lambda}} =  \sum_{\lambda} \frac{V(t_i : i \in \lambda)}{\prod_{k \in \lambda, l \notin \lambda}(t_l-t_k)},\]
where $V$ is a symmetric polynomial in $m$ variables. \newline

As in the case of projective space, we can write an expression for the righthand-side of the Berline-Vergne formula as in iterated residue, in the following way:
\[\int_{\Grass{m}{n}} \phi(\mathcal{R}) = Res_{z_1 = \infty}Res_{z_2 = \infty}...Res_{z_m = \infty} \frac{1}{m!}\frac{V(z_1,...,z_m)\prod_{i \neq j}(z_i - z_j)}{\prod_{i=1}^n\prod_{j=1}^m(t_i - z_j)}.\]
It can be shown directly, making similar calculations as for the projective space, and using the Residue Theorem $m$ times. We will proove a slightly more general result, from which the above formula follows easily.

\begin{formula} 
Let $\phi(\mathcal{R})$ be a characteristic class of the tautological bundle over $\Grass{m}{n}$, and $\psi(\mathcal{Q})$ a characteristic class of the quotient bundle. Then, the following fomula holds:

\[ \int_{\Grass{m}{n}} \phi(\mathcal{R})\psi(\mathcal{Q}) = \sum_{\lambda} \frac{\VR(t_k: k \in \lambda)\VQ(t_l: l \notin \lambda)}{\prod_{\substack{k \in \lambda, l \notin \lambda}} (t_l - t_k)}  =  \]
\[ =  Res_{\mathbf{z} = \infty}  \frac{ \VR(z_1,...,z_m)\VQ(z_{m+1},...,z_n)\prod_{i=1}^m\prod_{j=m+1}^n(z_i - z_j)\prod_{\substack{i,j =1 \\ i \neq j}}^m(z_i - z_j)\prod_{\substack{i,j =m+1 \\ i \neq j}}^n(z_i - z_j)}{m!(n-m)!\prod_{i=1}^{n}\prod_{j=1}^{n}(t_i - z_j)},\]
\label{classical}
\end{formula}
\begin{proof}
The first equality is just an application of Berline-Vergne formula. The lefthand-side can be rewritten using Vandermonde determinants. We have:
\[
\left\{  \begin{array}{l}\displaystyle\prod_{\substack{i,j =1 \\ i < j}}^m(z_i - z_j)\prod_{\substack{i,j =m+1 \\ i < j}}^n(z_i - z_j)\prod_{i=1}^m\prod_{j=m+1}^n(z_i - z_j) = (-1)^{mn}Vand(z_1,...,z_n)\\
\displaystyle\prod_{\substack{i,j =1 \\ i > j}}^m(z_i - z_j)\prod_{\substack{i,j =m+1 \\ i > j}}^n(z_i - z_j) = \frac{Vand(z_1,\cdots,z_n)}{\prod_{i=1}^m\prod_{j=m+1}^n(z_j - z_i)},\end{array} \right.\]
where $Vand(x_1,\cdots,x_n) =  \prod_{i<j}(x_i - x_j)$ denotes the Vendermonde determinant. \newline

Now, the lefthand-side can be written as:
\[ Res_{z_1,...,z_n=\infty} \frac{\VR(z_1,...,z_m)\VQ(z_{m+1},...,z_n)(-1)^{mn}Vand(z_1,...,z_n)^2(\prod_{i=1}^m\prod_{j=m+1}^n(z_j - z_i))^{-1}}{m!(n-m)!\prod_{i=1}^{n}\prod_{j=1}^{n}(t_i - z_j)} = \star\]
The residue at infinity can be calculated by replacing it with residues at points $t_i$, by the Residue Theorem. The numerator of $\star$ is a polynomial, so we only need to check what is the result of taking residue of the function $\frac{1}{\prod_{i=1}^{n}\prod_{j=1}^{n}(t_i - z_j)}$.To shorten the expressions appearing in the calculation, let's use the following notation:
\[ R_y^x(I,J):= \prod_{i \in I}\prod_{j \in J}(x_i - y_j) .\] 
In a special case, when $I=J$ and $x=y$, we will write $R_x(I)$ instead of $R_x^x(I,I)$.
We will also use the notation $[n]$ for $\{1,\dots,n\}$ and $[m,n]$ for $\{m, m+1,\dots, n\} $. \newline

Now, the residue of $\frac{1}{\prod_{i=1}^{n}\prod_{j=1}^{n}(t_i - z_j)}$ is:

\[Res_{z_1,...,z_n=\infty} \frac{1}{\prod_{i=1}^{n}\prod_{j=1}^{n}(t_i - z_j)} = Res_{z_1,...,z_n=\infty} \frac{1}{R_z^t([n],[n])} = \]
\[Res_{z_1,...,z_{n-1}=\infty} \sum_{k_n = 1}^n \frac{1}{R_z^t([n],[n-1])\prod_{i \neq k_n}(t_i - t_{k_n})} = \]
\[ =  Res_{z_1,...,z_{n-2}=\infty} \! \sum_{\substack{k_{n-1} = 1 \\ k_{n-1}\neq k_n}}^n \! \sum_{k_n = 1}^n  \frac{1}{ R_z^t([n],[n-2])R_t^t([n]\setminus\{k_n,k_{n-1}\},\{k_n,k_{n-1}\})R_t(\{k_n,k_{n-1}\})} = \]
\[ = \cdots = \sum_{\substack{k_1 = 1 \\ k_1 \neq k_2,...,k_n}}^n \cdots \sum_{k_n = 1}^n \frac{1}{R_z^t([n],\varnothing)R_t^t(\varnothing,\{k_1,\dots,k_n\})R_t(\{k_1\dots ,k_n\})}\]
\[=\sum_{\substack{k_1 = 1 \\ k_1 \neq k_2,...,k_n}}^n \cdots \sum_{k_n = 1}^n \frac{1}{\prod_{i,j = 1}^n (t_{k_i} - t_{k_j})} =  \sum_{\substack{k_1 = 1 \\ k_1 \neq k_2,...,k_n}}^n \! \cdots \sum_{k_n = 1}^n \frac{1}{(-1)^{mn}Vand(t_{k_1},...,t_{k_n})^2}\]

Finally, we obtain:
\[ \star = Res_{z_1,...,z_n=\infty} \frac{\VR(z_1,...,z_m)\VQ(z_{m+1},...,z_n)(-1)^{mn}Vand(z_1,...,z_n)^2\bigg(R_z^z([m+1,n],[m])\bigg)^{-1}}{m!(n-m)!\prod_{i=1}^{n}\prod_{j=1}^{n}(t_i - z_j)} = \]
\[ = \sum_{\substack{k_1 = 1 \\ k_1 \neq k_2,...,k_n}}^n \cdots \sum_{k_n = 1}^n \frac{\VR(t_{k_1},...,t_{k_m})\VQ(t_{k_{m+1}},...,t_{k_n})(-1)^{mn}Vand(t_{k_1},...,t_{k_n})^2}{m!(n-m)!(-1)^{mn}Vand(t_{k_1},...,t_{k_n})^2 (R_{t_k}^{t_k}([m+1,n],[m])} =  \]
\[ = \sum_{\substack{k_1 = 1 \\ k_1 \neq k_2,...,k_n}}^n \!\!\!\! \cdots \sum_{k_n = 1}^n \frac{\VR(t_{k_1},...,t_{k_m})\VQ(t_{k_{m+1}},...,t_{k_n})}{m!(n-m)!\prod_{i=1}^m\prod_{j=m+1}^n(t_{k_j} - t_{k_i})} =\]
\[= m! (n-m)!\cdot \sum_{\lambda} \frac{\VR(t_k: k \in \lambda)\VQ(t_l: l \notin \lambda)}{m!(n-m)!\prod_{\substack{k \in \lambda, l \notin \lambda}} (t_l - t_k)}, \]
where the last equality results from the fact, that $\VR$ and $\VQ$ are symmetric polynomials, so for any permutation $\sigma \in S_m \times S_{n-m}$ the residues at point $(t_{k_1}, \dots , t_{k_n})$ and $(t_{\sigma(k_1)}, \dots , t_{\sigma(k_n)})$ are equal. The final experession is exactly the righthand-side in the formulation of the theorem. 
\end{proof}

\begin{cor}
If we set $\VQ$ to be a constant polynomial, $\VQ=1$, we get
\[ \int_{\Grass{m}{n}} \phi(\mathcal{R}) = Res_{z_1,\dots,z_m = \infty} \frac{1}{m!}\frac{V(z_1,...,z_m)\prod_{i \neq j}(z_i - z_j)}{\prod_{i=1}^n\prod_{j=1}^m(t_i - z_j)}.\] 
\end{cor}
The proof of the above Corollary is a simple calculation. \newline

\subsection{Lagrangian Grassmannians}\label{lagrangian}
Similar results can be obtained for integrals over the Lagrangian Grassmannian. \newline

For $LG(n)$, the fixed points of the torus action can be parametized using the subsets $I \subseteq \{ 1,\dots,n\}$:
\[p_I = Span\{ q_i, p_j: i \in I, j \notin I \}.\]
Weights of the torus action on the tangent space are equal to
\[\{ \pm t_i \pm t_j: 1\leq i<j \leq n\}\cup\{\pm 2t_i: i=1,...,n\} \textrm{, where the + sign appears whenever } i,j \in I.\]
In this case, the Berline-Vergne formula gives:
\[ \int_{LG(n)} \phi(\mathcal{R}) = \sum_{I}\frac{V(t_i, -t_j: i \in I, j \notin I)}{\prod_{i,j=1}^n (\pm2t_i)( \pm t_i \pm t_j)}, \]
which can be expressed as a residue at infinity as follows:

\begin{formula}
\[\int_{LG(n)} \phi(\mathcal{R}) =  Res_{\mathbf{z} = \infty} \frac{V(z_1,...,z_n)\prod_{i<j}(z_j - z_i)}{\prod_{i=1}^n(t_i - z_i)(t_i + z_i)\prod_{i<j}(t_i + t_j)(t_j - t_i)}.\] \label{formula-lagrangian}
\end{formula}

\begin{proof}
First, we use the Residue Theorem to change the residue at infinity to residues at points $z_i = \pm t_i$. The function
\[\widetilde{V} = \frac{V(z_1,...,z_n)\prod_{i<j}(z_j - z_i)}{\prod_{i<j}(t_i + t_j)(t_j - t_i)}\] 
is a polynomial in $z_1,...,z_n$, so we only need to compute the residue of the function~$\frac{1}{\prod_{i=1}^n(t_i - z_i)(t_i + z_i)}$. We have:

\[  Res_{z_1 = \pm t_1} \cdots Res_{z_n = \pm t_n} \frac{1}{\prod_{i=1}^n(t_i - z_i)(t_i + z_i)} = \]
\[= Res_{z_1 = \pm t_1}\frac{1}{(t_1-z_1)(t_1+z_1)}\cdots Res_{z_n = \pm t_n}\frac{1}{(t_n-z_n)(t_n+z_n)} = \frac{1}{\pm2t_1}\cdots\frac{1}{\pm2t_n}.\]

Now it is sufficient to calculate the values $\widetilde{V}$ at points $z_i = t_i, i \in I $, $z_j = -t_j, j \notin I$:

\[\widetilde{V}(\pm t_1,...,\pm t_n) = \frac{V(\pm t_1,...,\pm t_n)\prod_{i<j}(\pm t_j - \pm t_i)}{\prod_{i<j}(t_i + t_j)(t_j - t_i)}. \]

The interpretation of the $\pm$ signs in the product $\prod_{i<j}(\pm t_j - \pm t_i)$ is the followinng:

\[\prod_{i<j}(\pm t_j - \pm t_i) = \prod_{\substack{i<j \\ i,j \in I}}(t_j - t_i) \prod_{\substack{i<j\\ i,j \notin I}}(-t_j + t_i) \prod_{\substack{i<j\\ i \in I, j \notin I}}(-t_j - t_i) \prod_{\substack{i<j\\ i \notin I, j \in I}}(t_j + t_i) \] 

Reducing the above product with the part $\prod_{i<j}(t_i + t_j)(t_j - t_i)$ in the denominator of~$\widetilde{V}$, we obtain that  $\widetilde{V}(t_i, -t_j: i \in I, j \notin I)$ is equal to:

\[ \frac{V(t_i, -t_j: i \in I, j \notin I)}{\prod_{\substack{i<j \\ i,j \in I}}(t_i + t_j) \prod_{\substack{i<j\\ i,j \notin I}}(t_i + t_j)(-1)^{C_1} \prod_{\substack{i<j\\ i \in I, j \notin I}}(t_j - t_i) \prod_{\substack{i<j\\ i \notin I, j \in I}}(t_j - t_i)(-1)^{C_2}}, \]

where $C_1 = \#\{ (i,j): i<j, i \notin I, j \in I\} $, $C_2 = \#\{ (i,j): i<j, i,j \in I\} $. \newline

Changing the signs of $t_j$ to negative for $j \notin I$ and reducing them with the $(-1)^{C_1+C_2}$, we can rewrite it in the following way:

\[\widetilde{V}(t_i, -t_j: i \in I, j \notin I) = \frac{V(t_i, -t_j: i \in I, j \notin I)}{\prod_{\substack{i<j \\ i,j \in I}}(t_i + t_j) \prod_{\substack{i<j\\ i,j \notin I}}(-t_i - t_j) \prod_{\substack{i<j\\ i \in I, j \notin I}}(-t_j + t_i) \prod_{\substack{i<j\\ i \notin I, j \in I}}(t_j - t_i)}. \]

Finally, 

\[Res_{\mathbf{z} = \infty} \frac{V(z_1,...,z_n)\prod_{i<j}(z_j - z_i)}{\prod_{i=1}^n(t_i - z_i)(t_i + z_i)\prod_{i<j}(t_i + t_j)(t_j - t_i)} = \sum_{I} \frac{\widetilde{V}(t_i, -t_j: i \in I, j \notin I)}{\pm2t_1...\pm2t_n} =  \]
\[ = \frac{V(t_i, -t_j: i \in I, j \notin I)}{\prod_{i \in I}2t_i\prod_{j \notin I}(-2t_j)\prod_{\substack{i<j \\ i,j \in I}}(t_i + t_j) \prod_{\substack{i<j\\ i,j \notin I}}(-t_i - t_j) \prod_{\substack{i<j\\ i \in I, j \notin I}}(-t_j + t_i) \prod_{\substack{i<j\\ i \notin I, j \in I}}(t_j - t_i)}, \]

which shows that the residue we calculated is equal to the sum from the Berline-Vergne formula.
\end{proof}

Formula \ref{formula-lagrangian} is the natural expression one comes up with, when trying to find the easiest way to write down  the integral over $LG(n)$ as a residue at infinity of some function. However, one can easily find a slightly different formula, which gives the same result, but shows more similarities to the formula for the classical Grassmannian:

\begin{formula}
\[\int_{LG(n)} \phi(\mathcal{R}) =   Res_{\mathbf{z} = \infty} \frac{1}{n!}\frac{V(z_1,...,z_n)\prod_{i \neq j}(z_j - z_i)\prod_{i < j}(z_i+z_j)}{\prod_{i=1}^n\prod_{j=1}^n(t_i - z_j)(t_i + z_j)}.\]
\end{formula}

This second version of the formula may not seem very different from the first one, and in no way a better one. But it turns out to be more useful - it will derive the general formula for all the Grasmannian types in section \ref{general}. For orthogonal Grasmannians we will use only the latter form. The proof is analogous to the one in the classical case.

\subsection{Orthogonal Grassmannians}\label{orthogonal} 
For orthogonal Grassmannians, the fixed points can be indexed by the subsets $I \subseteq \{ 1,\dots,n\}$, just like in the Lagrangian case. The fixed points are: 
\[p_I = Span\{ q_i, p_j: i \in I, j \notin I \}.\]
The weights of the torus action on the tangent space, and hence the equivariant Euler class, depend on whether we consider the even-dimensional case of $OG(n,2n)$ or the odd-dimensional $OG(n,2n+1)$. In the former case, the weights are:
\[\{ \pm t_i \pm t_j: 1\leq i<j \leq n\} \textrm{,  where the + sign appears whenever } i,j \in I.\]
Therefore, the Berline-Vergne formula yields
\[ \int_{OG(n,2n)} \phi(\mathcal{R}) = \sum_{I}\frac{V(t_i, -t_j: i \in I, j \notin I)}{\prod_{i,j=1}^n ( \pm t_i \pm t_j)}. \] 
The righthand-side of the above equality is almost identical to the one in Berline-Vergne formula for Lagrangian Grassmannian. In order to obtain an iterated residue formula for $OG(n,2n)$ we only need to modify slightly the result for the Lagrangian case. 

\begin{formula}
\[ \int_{OG(n,2n)} \phi(\mathcal{R}) =  Res_{\mathbf{z} = \infty} \frac{1}{n!}\frac{V(z_1,...,z_n)\prod_{i \neq j}(z_j - z_i)\prod_{i < j}(z_i+z_j)2^n\prod_{i=1}^n z_i}{\prod_{i=1}^n\prod_{j=1}^n(t_i - z_j)(t_i + z_j)}.\]
\end{formula}

As might be expceted, the odd-dimensional case is just as easy - the weights of the torus action are
\[\{ \pm t_i \pm t_j: 1\leq i<j \leq n\}\cup\{ t_i: i=1,...,n \} \textrm{, where the + sign appears whenever } i,j \in I.\]
Hence, the residue formula for $OG(n,2n+1)$ is the following:
\begin{formula}
\[ \int_{OG(n,2n+1)} \phi(\mathcal{R}) =  Res_{\mathbf{z} = \infty} \frac{1}{n!}\frac{V(z_1,...,z_n)\prod_{i \neq j}(z_j - z_i)\prod_{i < j}(z_i+z_j)2^n}{\prod_{i=1}^n\prod_{j=1}^n(t_i - z_j)(t_i + z_j)}.\]
\end{formula}

\subsection{General formula}\label{general}
Surprisingly, it turns out that for the established formulas, one can find one equality, which generalises them all. It requires using the Weyl group action on the characters of the natural representation of the torus action. \newline

The Grassmannians are homogenous spaces $G/P$, where $G$ is the general linear group $GL(n)$ for the classical Grassmannian, the symplectic group $Sp(\mathbb{C}^n)$ for the Lagrangian Grassmannian, and $OG(n)$ for the orthogonal Grassmannian. The parabolic subgroup $P$ in each case is the stabilizer of the point under the action of the torus $\mathbb{T}$ of diagonal matrices. Let $W$ denote the Weyl group of $G$, and $W_P$ the Weyl group of $P$. Let $t_i$ denote the characters of the natural representation of the torus, and consider the "second copy" of those characters, $z_i$, which we are going to treat like formal variables. \newline

Now, let $X_i$ denote the set of all possible images of $z_i$ under the action of $W$, namely $X_{i} = \{ \sigma(z_i): \sigma \in W \}$. Furthermore, let $\Phi^+\setminus{\Phi_P}^+$ be the set of positive roots, written in terms of characters $z_i$. The following formula holds for all the considered cases:

\begin{formula}
Let $G/P$ be the classical, Lagrangian or orthogonal Grassmannian and let $\mathcal{R}$ denote the tautological bundle, given by polynomial $V$ in Chern roots. Then, we have

\[ \int_{G/P} \phi(\mathcal{R}) =  Res_{\mathbf{z} = \infty} \frac{1}{|W_P|}\frac{V(z_1,...,z_n)\prod_{i=1}^n \prod_{x \in X_i \setminus \{z_i\} } (z_i - x)}{\prod_{i=1}^n\prod_{x \in X_i} (t_i - x ) \prod_{y \in \Phi^+ \setminus {\Phi_P}^+} y }.\]

\end{formula}

\begin{proof}
It's enough to compare the above formula with the ones obtained in previous chapters. The table below shows, what the sets $X_i$ and the positive roots $\Phi^+\setminus{\Phi_P}^+$ are in each case. 
\vspace{1em}

\centerline{
\begin{tabular}{|c|c|c|} \hline
$G/P$ & $X_i$ & $\Phi^+\setminus{\Phi_P}^+$ \\ \hline
$\mathbf{Grass_m(\mathbb{C}^n)}$ & $\{ z_1,...,z_n \}$ & $\{ z_i - z_j: 1\leq i \leq m < j \leq n \}$ \\
$\mathbf{LG(n)}$ & $\{ z_1,...,z_n,-z_1,...,-z_n  \}$ & $\{  z_i + z_j: 1\leq i < j \leq n \}\cup\{2z_i: 1 \leq i \leq n \}$ \\
$\mathbf{OG(n,2n)}$ & $\{ z_1,...,z_n,-z_1,...,-z_n  \}$ & $\{  z_i + z_j: 1\leq i < j \leq n \}$ \\
$\mathbf{OG(n, 2n+1)}$ & $\{ z_1,...,z_n,-z_1,...,-z_n  \}$ & $\{  z_i + z_j: 1\leq i < j \leq n \}\cup\{z_i: 1 \leq i \leq n \}$ \\
$\mathbf{Flag_m(\mathbb{C}^n)}$ & $\{ z_1,...,z_n\}$ & $\{ z_i - z_j: 1\leq i < j \leq m \}\cup$ \\
 & & $\{ z_i - z_j: 1 \leq i \leq m < j \leq n \}$\\
 \hline
\end{tabular}}
\end{proof}

The existence of a "general" formula, which works for all types of Grassmannians, and depends only on some algebraic objects associated with them (Weyl group, positive roots, characters of the natural representation) gives rise to a natural question: can we go further, and find a consistent way of expressing intergals over homogenous spaces (maybe under some additional assumpiotns) as residues at intinity of some functions? If yes, what is the meaning of those functions, do they have an interpretation? In the nearest future, I will investigate this problem further, starting with an attempt to obtain a similiar result for small exceptional groups, and then maybe generalising the results to other homogenous spaces~$G/P$.


\begin{thebibliography}{}
\bibitem[1]{puppe} C.Allday, V.Puppe \textit{Cohomological methods in transformation groups}, Cambridge Studies in Advanced Mathematics vol.32, Cambridge University Press, 1993
\bibitem[2]{ab} M.Atiyah, R.Bott \textit{The moment map and equivariant cohomology}, Topology vol.23, 1984, pp. 1-28
\bibitem[3]{berczi} G.B\'erczi, A.Szenes \textit{Thom polynomials of Morin singularities}, Annals of Mathematics Volume 175 (2012), Issue 2, pp.567-629
\bibitem[4]{ginzburg} V.Ginzburg \textit{Equivariant cohomology and Kähler geometry}, Functional Analysis and its applications vol. 21, 1987, pp. 19-34.
\bibitem[6]{quillen} D.Quillen \textit{The Spectrum of an Equivariant Cohomology Ring: I
}, The Annals of Mathematics, Second Series, Vol. 94, 1971, pp. 549-572
\bibitem[5]{weber} A.Weber \textit{Equivariant Chern classes and localization theorem}, Journal of Singularities Volume 5 (2012), 153-176
\end{thebibliography}
\end{document}